\documentclass{amsart}
\usepackage{amssymb}

\usepackage{graphicx}
\usepackage{amsmath}

\newtheorem{theorem}{Theorem}
\theoremstyle{plain}

\newtheorem{corollary}{Corollary}

\newtheorem{definition}{Definition}
\newtheorem{example}{Example}

\newtheorem{lemma}{Lemma}

\newtheorem{problem}{Problem}
\newtheorem{proposition}{Proposition}
\newtheorem{remark}{Remark}

\numberwithin{equation}{section}

\input{tcilatex}

\begin{document}
\author{P. Martins Rodrigues and J. Sousa Ramos}
\address{Department of Mathematics, \\
Instituto Superior T\'{e}cnico, Univ. Tec. Lisboa\\
E-mail: pmartins@math.ist.utl.pt, \\
sramos@math.ist.utl.pt}
\title[Bowen-Franks groups as conjugacy invariants]{Bowen-Franks groups as conjugacy invariants for $\mathbb{T}^{n}$
automorphisms}
\maketitle

\begin{abstract}
The role of generalized Bowen-Franks groups (BF-groups) as topological
conjugacy invariants for $\mathbb{T}^{n}$ automorphisms is studied.

Using algebraic number theory, a link is established between equality of
BF-groups for different automorphisms ($BF$-equivalence) and an identical
position in a finite lattice ($\mathcal{L}$-equivalence). Important cases of
equivalence of the two conditions are proved.

Finally, a topological interpretation of the classical BF-group $\mathbb{Z}%
^{n}/\mathbb{Z}^{n}(I-A)$ in this context is presented.
\end{abstract}

\bigskip\ 


{\small \bigskip }


\section{introduction}

Bowen-Franks groups arise as a simple but nevertheless important invariant
in symbolic dynamics. After its appearance in a paper by R. Bowen and J.
Franks \textbf{[Bw-Fr 77]}, the second author showed \textbf{[Fr 84]},
continuing the work of Parry and Sullivan \textbf{[Pa-Su 75]}, that it is a
complete invariant (apart from a sign) for flow equivalence of irreducible
subshifts of finite type. This was generalized to the reducible case by
Huang \textbf{[Hu 94]}, \textbf{[Hu 95]}. Recently, the subject was treated
in a difference perspective in \textbf{[Bo 02]} and \textbf{[Bo-Hu 03]}.
Flow equivalence will be considered in the last section of this paper.

In the particular case of subshifts induced -via a Markov partition- by
families of piecewise monotone interval maps, the Bowen-Franks groups
present rather strong regularities, which have been exploited for instance
in \textbf{[Al-Se-SR 96]}.

It is worth mentioning, although the subject won't be touched in this paper,
that another important theoretical setting for Bowen-Franks groups is the
one of K-theory of C*-algebras (see \textbf{[Cu-Kr 80]}).

In this note we want to address the role of Bowen-Franks groups as
invariants in another situation, close to symbolic dynamics: torus
automorphisms. The Bowen-Franks groups are, once again, invariants for
topological conjugacy and they describe in detail the algebraic structure of
periodic orbits.

In \textbf{[MR 96] }and \textbf{[MR-SR 99] }we presented some partial
results concerning exclusively the two-dimensional case, where very tight
recurrence relations are easily deduced. The generalization of these
properties to higher dimensions, as well as other aspects of the subject are
the object of current research and we will present new results in a
forthcoming paper.

Here, we bring into play algebraic number theory to show how good as an
invariant the Bowen-Franks groups are. In fact, Bowen-Franks groups reveal,
in a simple and clear manner, quite a lot of the arithmetic structure lying
behind these dynamical systems. We will show below that, at least for a
large class of automorphisms with the same characteristic polynomial, their
Bowen-Franks groups depend only on the position of the automorphism in a
certain finite lattice.

Just to illustrate the consequences, and referring the reader to the
definitions in the first section, we state the

\begin{proposition}
Let $p(x)\in \mathbb{Z}\left[ x\right] $ be an irreducible monic polynomial
of degree $n$, such that $p(0)=\pm 1$ and with square-free discriminant. Let 
$g(x)\in \mathbb{Z}\left[ x\right] $. Then all $\mathbb{T}^{n}$%
-automorphisms with characteristic polynomial $p(x)$ have the same
Bowen-Franks group $BF_{g}.$
\end{proposition}

This will be a corollary of our main result.

The results presented in this paper are related with symbolic dynamics since
generalized $BF$-groups are also, at least for polynomials satisfying $%
g(0)=\pm 1$, invariants of shift equivalence. In the particular case of
primitive, irreducible matrices with $\det =\pm 1$, shift equivalence
coincides with similarity and therefore the $BF$-equivalence classification
directly applies. We refer the reader to \textbf{[Li-Ma 95]} for all
definitions and results concerning symbolic dynamics.

The main interest of this work lies, from our point of view, on its
connection to the general topological classification of torus automorphisms.
Although there are algorithms that solve, in principle, this problem (%
\textbf{[Bo-Sh 66], [Gr-Se 80]}), the goal of finding a complete set of
computable and meaningful conjugacy invariants was attained, as far as we
know, only in dimension two (\textbf{[MR-SR 98], [Ad-Tr-Wo 97]}).

The links of torus automorphisms with algebraic number theory have also been
exploited in the context of $\mathbb{Z}^{d}$ actions (cf. \textbf{[Sc 95]}, [%
\textbf{Li-Sc 99]}, [\textbf{Ka-Ka-Sc 00}]).

\section{Basic facts}

We briefly recall some basic facts. The $n$-torus $\mathbb{T}^{n}=\mathbb{R}%
^{n}/\mathbb{Z}^{n}$ has $M_{n}(\mathbb{Z)}$ as its ring of endomorphisms
with group of units $GL_{n}(\mathbb{Z)}$. We will use the same notation for
a torus endomorphism and its associated matrix. Matrices will act on $%
\mathbb{T}^{n}$ on the left and so torus points will be represented by
column vectors. $\mathbb{Z}^{n}$ denotes the lattice of points with integer
coordinates, and are also represented by column vectors, while the dual
group of $\mathbb{T}^{n}$, whose elements are represented by row integer
vectors, will be denoted $\left( \mathbb{T}^{n}\right) ^{\ast }$.

Topological conjugacy of endomorphisms is equivalent to algebraic conjugacy,
as one easily deduces by considering the action of homeomorphisms on the
one-dimensional homology.

\begin{remark}
Although most of what we will say applies to general endomorphisms, we will
consider only irreducible automorphisms; this restriction is not important
for classification purposes (cf. the remark 16 below).
\end{remark}

The simplest algebraic invariant one can think of is, of course, the
characteristic polynomial. This determines the two basic dynamically
significant invariants, topological entropy and the zeta function, by
well-known results.

A first step beyond the mere counting of periodic points is given by means
of the following

\begin{definition}
Let $A$ be a square integer matrix of dimension $n$. The Bowen-Franks group
of $A$ is defined as $BF(A)=\mathbb{Z}^{n}/(A-I)\mathbb{Z}^{n}$.
\end{definition}

As a simple consequence of the structure theorem for finitely generated
abelian groups, we have $BF(A)=\mathbb{Z}_{a_{1}}\oplus \cdots \oplus 
\mathbb{Z}_{a_{m}}\oplus \mathbb{Z}^{n-m}$, with $a_{1}\mid a_{2}\mid \cdots
\mid a_{m}$, where the torsion coefficients $a_{i}$ and the Betti number $%
n-m $ may be computed through the reduction of $A-I$ to its Smith normal
form.

Bowen-Franks groups give rise not to one but to infinitely many conjugacy
invariants by generalizing the polynomial in the definition. Explicitly

\begin{definition}
Let $g(x)\in \mathbb{Q}\left[ x\right] $ such that $g(A)\in M_{n}(\mathbb{Z)}
$. The (generalized) Bowen-Franks group of $A$ associated to $g(x)$ is $%
BF_{g}(A)=\mathbb{Z}^{n}/g(A)\mathbb{Z}^{n}$.
\end{definition}

\begin{definition}
Two integer matrices $A$ and $B$ with the same characteristic polynomial are
called Bowen-Franks-equivalent (or just BF-equivalent) if $%
BF_{g}(A)=BF_{g}(B)$, $\forall g(x)\in \mathbb{Q}\left[ x\right] $ such that 
$g(A),g(B)\in M_{n}(\mathbb{Z)}$.
\end{definition}

\begin{remark}
It will become clear later that a necessary condition for $A$ and $B$ to be $%
BF$-equivalent is that $g(A)\in M_{n}(\mathbb{Z)\Leftrightarrow }g(B)\in
M_{n}(\mathbb{Z)}\forall g(x)\in \mathbb{Q}\left[ x\right] $.
\end{remark}

A standard argument shows the existence of an isomorphism of abelian groups 
\begin{equation*}
\mathbb{Z}^{n}/(A^{k}-I)\mathbb{Z}^{n}\approx Per_{k}(A)
\end{equation*}
where $Per_{k}(A)$ is the group of $k$-periodic points of $A$ in $\mathbb{T}%
^{n}$. Because of their importance we use the special notation $BF_{k}(A)$
for $\mathbb{Z}^{n}/(A^{k}-I)\mathbb{Z}^{n}$. This will cause no confusion
since $BF_{g}$ groups for constant polynomials are uninteresting and never
used.

\begin{example}
Let 
\begin{equation*}
A=\left( 
\begin{array}{ccc}
0 & 1 & 0 \\ 
0 & 0 & 1 \\ 
1 & -7 & 23
\end{array}
\right) ,B=\left( 
\begin{array}{ccc}
0 & 1 & 0 \\ 
1 & 0 & -3 \\ 
0 & 2 & 23
\end{array}
\right) ,C=\left( 
\begin{array}{ccc}
0 & 1 & 0 \\ 
1 & 0 & 4 \\ 
6 & -2 & 23
\end{array}
\right) \text{.}
\end{equation*}
These matrices have the same characteristic polynomial $%
p(x)=x^{3}-23x^{2}+7x-1$. Then $BF_{1}(A)=\mathbb{Z}_{16}$ while $%
BF_{1}(B)=BF_{1}(C)=\mathbb{Z}_{2}\oplus \mathbb{Z}_{8}$; but $BF_{x+1}(B)=%
\mathbb{Z}_{2}\oplus \mathbb{Z}_{16}$ and $BF_{x+1}(C)=\mathbb{Z}_{4}\oplus 
\mathbb{Z}_{8}$.
\end{example}

The torsion coefficients of $BF_{k}(A)$ enable us to locate the $k$-periodic
orbits: let $T_{m}=\left\{ x\in \mathbb{T}^{n}:mx=0\right\} $; we have

\begin{proposition}
Let $BF_{k}(A)=\mathbb{Z}_{k_{1}}\oplus \cdots \oplus \mathbb{Z}_{k_{n}}$
where $k_{1}\mid k_{2}\mid \cdots \mid k_{n}$ Then $Per_{k}(A)$ is generated
by $x_{i}\in Per_{k}(A)\cap T_{k_{i}}$, $i=1,\cdots ,n$; in particular $%
T_{k_{1}}\subset Per_{k}(A)\subset T_{k_{n}}$ and these are the best
possible inclusions.
\end{proposition}

\begin{remark}
In the previous proposition we are implicitly assuming that $\det
(A^{k}-I)\neq 0$, i.e. that $A$ has no $k$-roots of unity among its
eigenvalues. Of course the statement is still true in the general case if we
extend the definitions in the obvious way: $\mathbb{Z}_{0}\mathbb{=Z}$ and $%
T_{0}=\mathbb{T}^{n}$. In general, we will consider only Bowen-Franks groups 
$BF_{g}(A)$ such that $\det (g(A))\neq 0$.
\end{remark}

The above proposition is generalized to the case of an also generalized $%
BF_{g}(A)$ if we just replace $Per_{k}(A)$ by the group of periodic orbits 
\begin{equation*}
T_{g}=\left\{ x\in \mathbb{T}^{n}:g(A)x=0\right\} .
\end{equation*}

\section{BF-groups and algebraic numbers}

As we mentioned in the introduction, the study of $\mathbb{T}^{n}$%
-automorphisms benefits from the connection with Algebraic Number Theory. To
establish this connection in the simplest way, let's fix an irreducible
monic polynomial $p(x)\in \mathbb{Z}\left[ x\right] $. Consider the ring $%
\mathbb{Z}\left[ x\right] /(p(x))$ and its field of fractions $K=\mathbb{Q}%
\left[ x\right] /(p(x))$. For later use we recall the basic

\begin{definition}
Let $\theta \in K$ and $\theta _{1}=\theta ,\theta _{2},\ldots ,\theta _{n}$
\ be its conjugates. The absolute Trace and Norm of $\theta $ are defined
resp. as $Tr(\theta )=\sum_{i=1}^{n}\theta _{i}$ and $N(\theta
)=\prod_{i=1}^{n}\theta _{i}$.
\end{definition}

If the characteristic polynomial of $A\in M_{n}(\mathbb{Z)}$ is $p(x)$ then,
as a $K^{n}$-endomorphism, $A$ has the equivalence class of $x$ as a simple
eigenvalue (because we assume $p(x)$ to be irreducible). We will denote this
equivalence class from now on by $\beta $ and identify $\mathbb{Z}\left[ x%
\right] /(p(x))$ with $\mathbb{Z}\left[ \beta \right] $.

Let $v=\left( v_{1},\cdots ,v_{n}\right) $ be an associated row eigenvector;
its entries form a $\mathbb{Z}$-basis of a (possibly fractional) ideal $%
I\vartriangleleft \mathbb{Z}\left[ \beta \right] $. Recall that a fractional
ideal is just a finitely generated $\mathbb{Z}\left[ \beta \right] $-module
contained in $K$ and that these are exactly the sets $zJ$ for some $z\in K$
and $J\subset \mathbb{Z}\left[ \beta \right] $ an (integral) ideal. Two
ideals $I,J$ are said to be equivalent if $J=zI$ for some nonzero $z\in K$.
Thus, a different eigenvector $w=zv$ gives rise to an equivalent ideal. On
the other hand, passing from $A$ to a similar matrix amounts to perform a
change of base on $I$.

In short, we have just sketched the proof of the classic:

\begin{theorem}
(Latimer-McDuffee-Taussky) Given an irreducible polynomial $p(x)\in \mathbb{Z%
}\left[ x\right] $, there is a bijection between similarity classes of
integer matrices with characteristic polynomial $p(x)$ and ideal classes in $%
\mathbb{Z}\left[ x\right] /(p(x))$.
\end{theorem}

This theorem gives us the algebraic number theoretic version of the
classification problem (cf. [Nm 72] for a detailed proof).

With the previous notation, the mapping 
\begin{equation*}
\left( x_{1},\cdots ,x_{n}\right) ^{t}\rightarrow \sum v_{i}x_{i}
\end{equation*}
yields bijections between $\mathbb{Q}^{n}$ and $K$, $\mathbb{Z}^{n}$ and $I$%
, and between $\mathbb{Q}^{n}/\mathbb{Z}^{n}$ and the torsion $\mathbb{Z}%
\left[ \beta \right] $-module $K/I$, where the $R$-module structure of $%
\mathbb{Q}^{n}$ is given by $(\alpha ,q)\in R\times \mathbb{Q}%
^{n}\rightarrow \alpha (A)q$.

A more geometrical approach involves considering all the real and complex
embeddings of $K$ (cf., for instance, \textbf{[Ke-Ve 98]}) but this will be
enough for our purposes.

The matrix $A$ represents multiplication by $\beta $ in $K$, and in
particular in $I$, with respect to the basis $v$. When there is no risk of
confusion, we'll speak simply of the associated matrix (or automorphism) of
an ideal or of the associated ideal of a matrix (or automorphism). We
stress, however, that this identification depends on the choice of $v$.

Now we have a new, obvious, interpretation for Bowen-Franks groups of $A$: 
\begin{equation*}
BF_{g}(A)\approx I/g(\beta )I
\end{equation*}

We will use these identifications without further explanation in the sequel.
Besides, as the polynomials $g(x)$ are taken $\func{mod}p(x)$, and may be
thus identified with elements of $K$, we will also use the notation $%
BF_{\alpha }(A)$ if $\alpha =g(\beta )$.

In particular, proposition 2 translates into the optimal inclusions 
\begin{equation*}
k_{1}I\supset g(\beta )I\supset k_{n}I
\end{equation*}

An immediate observation is that for every automorphism $A$ (equivalently,
for every $\mathbb{Z}\left[ \beta \right] $-ideal $I$) 
\begin{equation*}
\left| BF_{\alpha }(A)\right| =\left| I/\alpha I\right| =\left| N(\alpha
)\right| \text{.}
\end{equation*}

It's worthwhile to derive from this that $I/\alpha I\approx I/\lambda
I\Longrightarrow \left| \frac{N(\alpha )}{N(\lambda )}\right| =1$ and so if $%
\frac{\alpha }{\lambda }\in R$, we have in fact $\frac{\alpha }{\lambda }\in
R^{\ast }$. Recip., if $\eta =\frac{\alpha }{\lambda }\in R^{\ast }$, the
induced mapping 
\begin{eqnarray*}
\overset{\symbol{126}}{\eta } &:&I/\lambda I\rightarrow I/\alpha I \\
x+\lambda I &\rightarrow &\eta x+\alpha I
\end{eqnarray*}
is an isomorphism. In particular, we conclude that for a given automorphism $%
A$ only one group of each possible cardinality may occur as $BF_{\alpha }(A)$%
.

\subsection{Strong $BF$-equivalence}

In the definition of $BF$-equivalence of two automorphisms $A$ and $B$ given
before, $BF_{g}(A)$ and $BF_{g}(B)$ are demanded to be isomorphic only as
abelian groups. Another equivalence relation is

\begin{definition}
Two automorphisms $A$ and $B$ are strongly $BF$-equivalent if, for every
polynomial $g(x)$, $BF_{g}(A)$ and $BF_{g}(B)$ are isomorphic $\mathbb{Z}%
\left[ \beta \right] $-modules.
\end{definition}

In fact $BF_{g}(A)\simeq T_{g}(A)$ has the $\mathbb{Z}\left[ \beta \right] $%
-module structure induced by the action of $A$ making this a natural
definition: it requests that the isomorphism between these finite groups of
periodic orbits of $A$ and $B$ conjugates the actions of the automorphisms.

We will use the notation $BF_{g}(A)\simeq _{R}BF_{g}(B)$ to indicate the
existence of such an isomorphism and $A\simeq _{R}B$ to denote the strong $%
BF $-equivalence of the two automorphisms.

If $BF_{g}(A)\simeq _{R}BF_{g}(B)$ then the restriction of the isomorphism
to any submodule is again an isomorphism; if this submodule is of the form $%
T_{h}(A)=\left\{ x\in \mathbb{T}^{n}:h(A)x=0\right\} $ we conclude that $%
BF_{h}(A)\simeq _{R}BF_{h}(B)$. This simple observation shows that we may
have $BF_{g}(A)\simeq BF_{g}(B)$ and not $BF_{g}(A)\simeq _{R}BF_{g}(B)$:

\begin{example}
Consider the two automorphims 
\begin{equation*}
M=\left( 
\begin{array}{rrrr}
-1 & -1 & -1 & -4 \\ 
4 & 1 & 3 & 8 \\ 
0 & 1 & 0 & 0 \\ 
0 & 0 & 1 & 7
\end{array}
\right) \text{ and  }M^{\prime }=\left( 
\begin{array}{rrrr}
-5 & -4 & -5 & -12 \\ 
8 & 5 & 7 & 12 \\ 
0 & 1 & 0 & 0 \\ 
0 & 0 & 1 & 7
\end{array}
\right) 
\end{equation*}
with characteristic polynomial $p(x)=$\ $x^{4}-7x^{3}-7x+1$; one verifies
that 
\begin{eqnarray*}
BF_{48}(M) &=&BF_{48}(M^{\prime })=\mathbb{Z}_{448}\oplus \mathbb{Z}%
_{1344}\oplus \mathbb{Z}_{m}\oplus \mathbb{Z}_{m} \\
m &=&130401445122840192
\end{eqnarray*}
Let's suppose that $BF_{48}(M)\simeq _{R}BF_{48}(M^{\prime })$, that is,
there is an isomorphism 
\begin{equation*}
\Psi :T_{h}(M)=\left\{ x\in \mathbb{T}^{4}:h(M)x=0\right\} \rightarrow
T_{h}(M^{\prime })=\left\{ x\in \mathbb{T}^{4}:h(M^{\prime })x=0\right\} 
\end{equation*}
where $h(x)=x^{48}-1$, such that $M^{\prime }\circ \Psi =\Psi \circ M$.

For $g(x)=x^{3}+4x^{2}+4x+5$, we have $BF_{g}(M)=\mathbb{Z}_{4}\oplus 
\mathbb{Z}_{8}\oplus \mathbb{Z}_{8}\oplus \mathbb{Z}_{64}$; in particular,
by proposition 2, 
\begin{equation*}
T_{g}(M)\subset T_{64}=\left\{ x\in \mathbb{T}^{4}:64x=0\right\} \subset
T_{h}(M)
\end{equation*}
and so $\Psi $ would induce an isomorphism between $T_{g}(M)$ and $%
T_{g}(M^{\prime })$.

However $BF_{g}(M^{\prime })=\mathbb{Z}_{8}\oplus \mathbb{Z}_{8}\oplus 
\mathbb{Z}_{8}\oplus \mathbb{Z}_{32}\neq BF_{g}(M)$.
\end{example}

This argument proves in fact the following

\begin{lemma}
Given automorphisms $A$ and $B$, if $BF_{k}(A)\simeq _{R}BF_{k}(B)$ $\forall
k$, then $A\simeq _{R}B$.
\end{lemma}

Observe that, in the previous example, the hypothesis 
\begin{equation*}
BF_{48}(M)\simeq _{R}BF_{48}(M^{\prime })
\end{equation*}
would imply that the restrictions of the actions by $M$ and $M^{\prime }$ to 
$T_{64}$ would be conjugate; this is the same as to say that the matrices $M$
and $M^{\prime }$ are conjugate in $GL_{4}(\mathbb{Z}_{64})$. On the other
hand, conjugacy $\func{mod}64$ would imply $BF_{g}(M)\simeq
_{R}BF_{g}(M^{\prime })$.

Since, given any automorphism $A\in GL_{n}(\mathbb{Z})$ and any natural
number $m$, there exists some $k$ such that every point in $T_{m}$ is $k$%
-periodic, the fact that for all $g(x)$, 
\begin{equation*}
BF_{g}(A)=\oplus _{i=1}^{n}\mathbb{Z}_{k_{i}}\Rightarrow T_{g}(A)\subset
T_{k_{n}}
\end{equation*}
and the previous lemma allow us to conclude that

\begin{proposition}
Two automorphisms $A$ and $B$ (with the same characteristic polynomial) are
strongly $BF$-equivalent iff they are conjugate $\func{mod}m$ $\forall m\in 
\mathbb{Z}^{+}$.
\end{proposition}

So, as we will see in the next sections, it can happen that two
automorphisms are conjugate in $GL_{n}(\mathbb{Z}_{m})$ $\forall m$ but not
in $GL_{n}(\mathbb{Z})$. In connection with this observation we quote the
following result (cf. \textbf{[Su 65]}):

\begin{theorem}
(D.A. Suprunenko) Let $A$ and $B$ be two square matrices of the same
dimension, over a unique factorization domain $R$ with a infinity of non
associate primes $p$. If $A$ and $B$ are similar over $R/(p)$ for all $p$,
then they are similar over the field of fractions of $R$.
\end{theorem}

\subsection{The lattice of orders and $\mathcal{L}-$equivalence}

Some more information from algebraic number theory will be needed. We'll
present it along with the development of our results. Although we keep this
to a minimum and refer the reader to the many textbooks on the subject (cf.
for instance, \textbf{[Co 78], [Fr-Ta 93], [Nk 99]}), most expositions on
algebraic number theory overlook the ideal theory of non-integrally closed
rings of integers (\textbf{[Nk 99]} is the best exception that we know
about) and so some useful results are either missing in the literature or
scattered, sometimes ``in disguise'', in treatises on commutative algebra as 
\textbf{[Ja 80]}. Thus, it seems adequate to include a detailed presentation
of these.

Recall that $z\in K$ is an integer of the field if its minimal polynomial
has rational integer coefficients. The set of all integers in $K$ is a ring,
denoted from now on as $Z_{K}$. This ring is known to be a Dedekind domain,
i.e.

i) all non zero prime ideals are maximal;

ii) the ascending chain condition for ideals is valid;

iii)$Z_{K}$ is integrally closed in $K$.

The main consequence of these properties for algebraic number theory is the
well known fact that any nonzero ideal of $Z_{K}$ is uniquely factorable as
a product of prime ideals. In our interest, let's just point out that all $%
Z_{K}$-ideals are invertible. Concretely, if $I$ is a fractional ideal,
define $I^{-1}=\left\{ z\in K:zI\subset Z_{K}\right\} $ which is again a
fractional ideal; then $II^{-1}=Z_{K}$.

The ring $\mathbb{Z}\left[ \beta \right] $ will not be, in general,
integrally closed in its field of fractions $K$. In other words, $\mathbb{Z}%
\left[ \beta \right] \subset Z_{K}$ will be a strict inclusion. In fact, $%
\mathbb{Z}\left[ \beta \right] $ and $Z_{K}$ are resp. the minimal and
maximal elements in a finite lattice, ordered by inclusion, of \textit{orders%
}, i.e. sub-rings of $Z_{K}$ \ with finite index and containing 1. Given $%
\mathbb{Z}\left[ \beta \right] $, or, more generally, some order, one may
construct explicitly its lattice (c.f. \textbf{[Co 78]}, for a description
of the algorithm).

\begin{example}
Let $p(x)=x^{2}-34x+1$; the discriminant is $(24)^{2}2$; the lattice
associated to $\mathbb{Z}\left[ \beta \right] $ is 
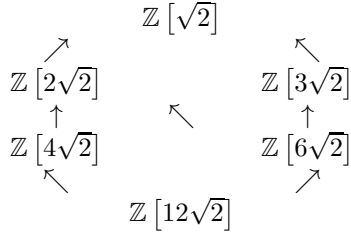
\begin{equation*}
\begin{array}{ccc}
& \mathbb{Z}\left[ \sqrt{2}\right] &  \\ 
\nearrow &  & \nwarrow \\ 
\mathbb{Z}\left[ 2\sqrt{2}\right] &  & \mathbb{Z}\left[ 3\sqrt{2}\right] \\ 
\uparrow & \nwarrow & \uparrow \\ 
\mathbb{Z}\left[ 4\sqrt{2}\right] &  & \mathbb{Z}\left[ 6\sqrt{2}\right] \\ 
\nwarrow &  & \nearrow \\ 
& \mathbb{Z}\left[ 12\sqrt{2}\right] & 
\end{array}
\end{equation*}
where the arrows describe inclusion of orders.
\end{example}

\begin{example}
Let $p(x)=x^{3}-23x^{2}+7x-1$; the discriminant is $-2^{8}83$; we have in
this case 
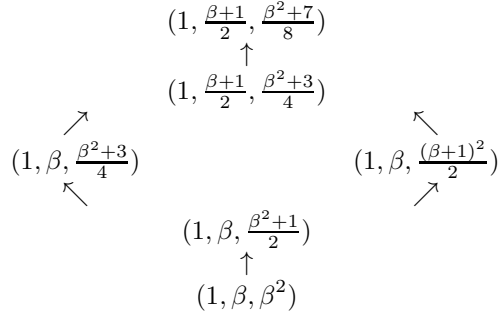
\begin{equation*}
\begin{array}{ccc}
& (1,\frac{\beta +1}{2},\frac{\beta ^{2}+7}{8}) &  \\ 
& \uparrow &  \\ 
& (1,\frac{\beta +1}{2},\frac{\beta ^{2}+3}{4}) &  \\ 
\nearrow &  & \nwarrow \\ 
(1,\beta ,\frac{\beta ^{2}+3}{4}) &  & (1,\beta ,\frac{(\beta +1)^{2}}{2})
\\ 
\nwarrow &  & \nearrow \\ 
& (1,\beta ,\frac{\beta ^{2}+1}{2}) &  \\ 
& \uparrow &  \\ 
& (1,\beta ,\beta ^{2}) & 
\end{array}
\end{equation*}
In this case, an integral basis of the corresponding order is displayed in
each vertex.
\end{example}

\begin{remark}
Notice that in the quadratic case the lattice is completely determined by
the discriminant: if $disc(p(x))=f^{2}\Delta $ where $\Delta $ is
square-free, then the lattice will be the one of the divisors of $f$ if $%
\Delta =1\func{mod}4$, and of $f/2$ otherwise. This doesn't happen for
higher dimensions.
\end{remark}

This lattice, denoted as $\mathcal{L}$ in the sequel, will be the
fundamental structure to identify BF-equivalence classes and we will pay
considerable attention to it.

It will be convenient to state some properties that hold for ideals in any
order $R\subset Z_{K}$, or, in certain cases, even in any noetherian domain.

Let $I$ be an ideal in $R$ and define, just as before, $I^{-1}=\left\{ z\in
K:zI\subset R\right\} $; we say that $I$ is invertible in $R$ if $II^{-1}=R$%
. Also, $I$ is said to be divisorial if $(I^{-1})^{-1}=I$. Notice that these
definitions depend on the ring in question. So, when necessary, we might use
the notation $(M:N)=\{z\in K:zN\subset M\}$ valid for any $R$-modules, where 
$R$ is any ring with field of fractions $K$.

The ring of coefficients of $I$ is defined as $C(I)=\left\{ z\in K:zI\subset
I\right\} $. By the Cayley-Hamilton theorem, $C(I)$ is contained in $Z_{K}$
and it obviously contains $R$, i.e. it's one of the orders in the lattice
described above.

If $I$ is the ideal associated to an automorphism (represented by the
matrix) $A$ one easily concludes that

\begin{lemma}
$R\subset C(I)$ iff $g(A)\in M_{n}[\mathbb{Z]}$, for every $g\in \mathbb{Q}%
[x]$ such that $g(\beta )\in R$.
\end{lemma}

This verification amounts to a trivial computation using the elements in the
basis of the different orders above $\mathbb{Z}\left[ \beta \right] $.

\begin{example}
Consider again the automorphisms $A,B,C$ from example 1; the rings of
coefficients of the associated ideals are the orders with integral basis $%
(1,\beta ,\beta ^{2}),$ $(1,\beta ,\frac{\beta ^{2}+1}{2})$ and $(1,\beta ,%
\frac{\beta ^{2}+3}{4})$\ resp.
\end{example}

\begin{definition}
Two automorphisms $A,A^{\prime }$ with the same characteristic polynomial
are $\mathcal{L}$-equivalent if their associated (classes of) ideals have
the same ring of coefficients.
\end{definition}

\begin{remark}
Before we return to BF-groups, we take the opportunity to show how the
classification of irreducible endomorphisms may be reduced to the one of
automorphisms. By the Latimer-McDuffee-Taussky theorem, the similarity
classes of endomorphisms with characteristic polynomial $q(x)$ are in
bijection with classes of ideals in $\mathbb{Z}\left[ x\right] /(q(x))$.

If $\beta $ is a unit of this ring, $\mathbb{Z}\left[ x\right] /(q(x))$ is
one of the orders above $\mathbb{Z}\left[ \beta \right] $ and its classes of
ideals are also classes of ideals in $\mathbb{Z}\left[ \beta \right] $, as
an equivalence of ideals in the former ring may always be turned into an
equivalence in the later. So, in order to identify the ideal classes of the
larger ring, one just has to choose among the ideal classes of $\mathbb{Z}%
\left[ \beta \right] $ those with coefficient ring containing $\mathbb{Z}%
\left[ x\right] /(q(x))$.
\end{remark}

\begin{example}
Let $q(x)=x^{2}-10x-7$; a fundamental unit in $\mathbb{Z}\left[ x\right]
/(q(x))=\mathbb{Z}\left[ 4\sqrt{2}\right] $ is $\beta =17+12\sqrt{2}$ whose
characteristic polynomial is $p(x)=x^{2}-34x+1$ as in example 3.

In dimension two, a complete and effective classification is known; from the 
$12$ ideal classes of $\mathbb{Z}\left[ \beta \right] $, one finds that only 
$3$ of them have coefficient ring containing $\mathbb{Z}\left[ 4\sqrt{2}%
\right] $. The same algorithm enable us to list canonical representatives of
the different conjugacy classes of automorphisms with characteristic
polynomial $p(x)$; to get a similar result for the ring of endomorphisms in
question, we have only to notice that these endomorphisms represent
multiplication by $5+4\sqrt{2}=\frac{\beta -2}{3}$ in their ideal classes;
if $A$ is an automorphism with characteristic polynomial $p(x)$, associated
to one of the $3$ ideal classes mentioned above, $\frac{A-2I}{3}$ is an
endomorphism with characteristic polynomial $q(x)$ associated to the same
class.
\end{example}

\section{$BF$-equivalence and $\mathcal{L}$-equivalence}

\subsection{First results}

Suppose that, for two different $\mathbb{Z}\left[ \beta \right] $-ideals $I$
and $J$, we have $C(I)\neq C(J)$. We may take, in that case some $\frac{%
\alpha }{m}\in C(I)\backslash C(J)$ with $\alpha \in R$ and $m\in \mathbb{Z}$%
. Then $\alpha I\subset mI$ and so $m\mid k_{1}$ if $I/\alpha I\approx 
\overset{n}{\underset{i=1}{\oplus }}\mathbb{Z}_{k_{i}}$. On the other hand, $%
\exists y\in J:\frac{\alpha }{m}y\notin J$; let $J/\alpha J\approx \overset{n%
}{\underset{i=1}{\oplus }}\mathbb{Z}_{l_{i}}$; we have the inclusion $\alpha
J\subset l_{1}J$, ie.. 
\begin{equation*}
\forall z\in J,\quad \exists w_{z}\in J:\alpha z=l_{1}w_{z}
\end{equation*}
. If $m\mid l_{1}$, then $\frac{\alpha }{m}y=\frac{l_{1}}{m}w_{y}\in J$, a
contradiction. We conclude that

\begin{lemma}
If $C(I)\neq C(J)$, then $BF_{\alpha }(I)\neq BF_{\alpha }(J)$ for some $%
\alpha \in \mathbb{Z}\left[ \beta \right] $.
\end{lemma}

Or, in terms of the definitions introduced before,

\begin{proposition}
If two automorphisms $A,A^{\prime }$ with the same characteristic polynomial
are BF-equivalent then they are $\mathcal{L}$-equivalent.
\end{proposition}

A similar reasoning leads us to a partial converse:

\begin{lemma}
If $C(I)=R$ and, for $\alpha \in R$, the Bowen-Franks groups are $BF_{\alpha
}(I)=\overset{n}{\underset{i=1}{\oplus }}\mathbb{Z}_{k_{i}}$ and $BF_{\alpha
}(R)=\overset{n}{\underset{i=1}{\oplus }}\mathbb{Z}_{l_{i}}$, then $%
k_{1}=l_{1}$, $k_{n}=l_{n}$.
\end{lemma}

\begin{proof}
$\alpha R\subset l_{1}R\Rightarrow \alpha I\subset l_{1}I$, therefore $%
l_{1}\leq k_{1}$; but $\alpha I\subset k_{1}I\Leftrightarrow \frac{\alpha }{%
k_{1}}\in C(I)=R\Rightarrow \alpha R\subset k_{1}R$ and so $k_{1}\leq l_{1}$%
. To conclude for the equality of the last coefficients we just have to
reverse inclusions.
\end{proof}

For (very) low $n$ this, together with the observation already made about
the order of BF-groups, immediately implies equivalence:

\begin{corollary}
Let $n\leq 3$. Two $\mathbb{T}^{n}$-automorphisms $A,A^{\prime }$ with the
same characteristic polynomial are $BF$-equivalent iff they are $\mathcal{L}$%
-equivalent.
\end{corollary}

We see that to arrive to a general result of this type, a detailed
description of the $\mathbb{Z}\left[ \beta \right] $-classes of ideals with
different coefficient rings is needed. As it will be seen later, the case of
invertible ideals is easier to handle.

\subsection{Invertible ideals}

It's a consequence of the definitions that if $I$ is invertible in $R$, then 
$C(I)=R$. The converse doesn't hold in general, as example 7 will show.

We start with two lemmas, the first of which may be found, for instance, in 
\textbf{[Ja 80].}

\begin{lemma}
If $I\vartriangleleft R$ is a proper ideal then $I^{-1}$ strictly contains $%
R $.
\end{lemma}

\begin{lemma}
Given $I\vartriangleleft R$ one has $%
C(I^{-1})=C((I^{-1})^{-1})=(II^{-1})^{-1}$.
\end{lemma}

\begin{proof}
1) Let $x\in C(I^{-1})$ i.e. $xI^{-1}\subset I^{-1}$. If $z\in \left(
I^{-1}\right) ^{-1}$ then $\forall y\in I^{-1},xzy\in zI^{-1}\subset R$.
Therefore $xz\in \left( I^{-1}\right) ^{-1}$ and so $x\in C(\left(
I^{-1}\right) ^{-1}).$

2) Now let $x\in C(\left( I^{-1}\right) ^{-1})$; given $a\in I$ and $y\in
I^{-1}$, since by definition $I\subset \left( I^{-1}\right) ^{-1}$ we have $%
xa\in \left( I^{-1}\right) ^{-1}$ and $xay\in I^{-1}\left( I^{-1}\right)
^{-1}\subset R$; this means $xII^{-1}\subset R$ , that is $x\in
(II^{-1})^{-1}$. \ 

3) Finally, let $x\in (II^{-1})^{-1}$; for $y\in I^{-1}$and $\forall a\in I$%
, $xya\in R$, i.e. $xy\in I^{-1}$ and $x\in C(I^{-1})$.
\end{proof}

This yields a characterization of invertibility for ideals in an order $R$
(which, in fact, is valid for any noetherian domain):

\begin{proposition}
An ideal $I\vartriangleleft R$ is invertible iff $C(I)=R$ and $I$ is
divisorial, iff $C(I^{-1})=R$.
\end{proposition}

We show now that in general there may exist $\mathbb{Z}\left[ \beta \right] $%
-ideals which are not invertible in its coefficient ring:

\begin{example}
Let $p(x)=x^{3}-23x^{2}+7x-1$; $(1,\beta ,\frac{\beta ^{2}+1}{2})$ is an
integral basis of $R$, one of the orders above $\mathbb{Z}\left[ \mathbb{%
\beta }\right] =\mathbb{Z}\left[ x\right] /(p(x))$, and the matrix
representing multiplication by $\beta $ with respect to it is 
\begin{equation*}
B=\left( 
\begin{array}{ccc}
0 & -1 & -11 \\ 
1 & 0 & -3 \\ 
0 & 2 & 23
\end{array}
\right) ;
\end{equation*}
Consider two other matrices with the same characteristic polynomial 
\begin{equation*}
A=\left( 
\begin{array}{ccc}
-7 & -7 & -20 \\ 
8 & 7 & 0 \\ 
0 & 1 & 23
\end{array}
\right) \text{ and }D=\left( 
\begin{array}{ccc}
-1 & -1 & -8 \\ 
2 & 1 & -6 \\ 
0 & 1 & 23
\end{array}
\right) \text{;}
\end{equation*}
the associated ideals $I$ and $J$ have integral basis $(8,\beta +7,\beta
^{2}+7)$ and $(2,\beta +1,\beta ^{2}+1)$ resp. Both have $R$ as ring of
coefficients. The inverses are $I^{-1}$ with basis $(1,\frac{\beta +1}{2},%
\frac{\beta ^{2}+2\beta +9}{16})$ and $J^{-1}$ with basis $(1,\frac{\beta -1%
}{2},\frac{(\beta +1)^{2}}{4})$; $C(I^{-1})=R$ but $C(J^{-1})=J^{-1}$ which
is a bigger order. Therefore $I$ is invertible but $J$ is not.
\end{example}

A somewhat different approach to invertibility rests on the fundamental
concept of localization. Recall that given a prime $P\vartriangleleft R$ the
localization of $R$ at $P$ is the subring of $K$, $R_{P}=\left\{ \frac{x}{y}%
:x\in R,y\in R\backslash P\right\} $. $R_{P}$ has a single maximal ideal $%
PR_{P}$. We quote two important results from \textbf{[Nk 99].}

\begin{proposition}
If $J\neq 0$ is an ideal of $R$, then 
\begin{equation*}
R/J\approx \underset{P}{\oplus }R_{P}/JR_{P}=\underset{P\supseteq J}{\oplus }%
R_{P}/JR_{P}\text{.}
\end{equation*}
\end{proposition}

\begin{proposition}
A fractional ideal $I$ of $R$ is invertible iff, for every prime $P\neq 0$, 
\begin{equation*}
I_{P}=IR_{P}
\end{equation*}
is a fractional principal ideal of $R_{P}$.
\end{proposition}

The following proposition is usually stated only for Dedekind domains. The
proof we present here adapts to the general case the one suggested in 
\textbf{[Nk 99].}

\begin{proposition}
If $I\vartriangleleft R$ is an invertible ideal, then $I$ has rank $2$ as a $%
R$-module.
\end{proposition}

\begin{proof}
Let $s\in I$; It suffices to show that $I/(s)$ is a principal $R/(s)$ ideal;
in fact, if $I/(s)=tR/(s)$ then $I=sR+tR$. But, taking $J=(s)$ in
proposition 6\ we have that $I/(s)$ is isomorphic to $\underset{P\supseteq
(s)}{\oplus }I_{P}/sR_{P}$ and by proposition 7\ , each $I_{P}$ is a
principal $R_{P}$-ideal. Therefore, $I/(s)$ is principal.
\end{proof}

\subsection{$BF$-equivalence: invertible ideals}

Recall the interpretation of Bowen-Franks groups made before: $%
BF_{g}(A)\simeq I/g(\beta )I$.

Given $\alpha ,\lambda \in C(I)=R$ there is an exact sequence 
\begin{equation*}
0\rightarrow I/\alpha I\rightarrow I/\alpha \lambda I\rightarrow I/\lambda
I\rightarrow 0
\end{equation*}
where the injection is $x+\alpha I\rightarrow \lambda x+\alpha \lambda I$
and the surjection is the canonical projection.

From the commutative diagram 
\begin{equation*}
\begin{array}{ccccccccc}
&  & 0 &  & 0 &  & ((\alpha )\cap I)/\alpha I &  &  \\ 
&  & \downarrow &  & \downarrow &  & \downarrow &  &  \\ 
0 & \rightarrow & I & \overset{\alpha }{\rightarrow } & I & \rightarrow & 
I/\alpha I & \rightarrow & 0 \\ 
&  & \downarrow &  & \downarrow &  & \downarrow &  &  \\ 
0 & \rightarrow & R & \overset{\alpha }{\rightarrow } & R & \rightarrow & 
R/\alpha R & \rightarrow & 0 \\ 
&  & \downarrow &  & \downarrow &  & \downarrow &  &  \\ 
&  & R/I & \overset{\alpha }{\rightarrow } & R/I &  & G &  & 
\end{array}
\end{equation*}
and applying the ``snake lemma'' we get 
\begin{equation*}
\begin{array}{ccccccccccc}
&  &  & \nearrow & I/\alpha I & \rightarrow & R/\alpha R & \searrow &  &  & 
\\ 
0 & \rightarrow & ((\alpha )\cap I)/\alpha I &  &  &  &  &  & G & \rightarrow
& 0 \\ 
&  &  & \searrow & R/I & \rightarrow & R/I & \nearrow &  &  & 
\end{array}
\end{equation*}
which tells us that in order to have a $R$-isomorphism $I/I\eta \approx
R/\eta R$ it suffices that $\eta $ induces a automorphism of $R/I$ (in other
words, that $\eta $ is a $R$-unity $\func{mod}I$).

Let's suppose that $I$ is invertible and $\alpha \in R$. The proof of
proposition 35\ shows we may take $\alpha $ as one of the generators of $%
I^{-1}$. If $\gamma $ is the other one, there exist $x,y\in I$ such that $%
\gamma x+\alpha y=1$; the ideal $J=\gamma I$ is also contained in $R$ and
satisfies $J+\alpha R=R$ implying that $J\cap \alpha R=\alpha J$. In
consequence the morphism $J/\alpha J\rightarrow R/\alpha R$ is injective
and, as the modules are finite, bijective.

Since $I/\alpha I\simeq J/\alpha J$ we proved the following:

\begin{proposition}
If $A$ and $B$ are automorphisms associated resp. to the invertible ideal $I$
and to its coefficient ring $R$, then $A$ and $B$ are strongly $BF$%
-equivalent.
\end{proposition}

\begin{remark}
The above proof shows that $BF_{\alpha }(A)\simeq _{R}BF_{\alpha }(B)$, if
there exists $\gamma \in I^{-1}$ such that $\gamma I+\alpha R=R$.
\end{remark}

If $\mathbb{Z}\left[ x\right] /(p(x))=\mathbb{Z}\left[ \mathbb{\beta }\right]
$ is the ring $Z_{K}$ of all algebraic integers of $K$, we obtain that all
automorphisms with characteristic polynomial $p(x)$ are strongly $BF$%
-equivalent. That's what happens, in particular, if the discriminant of $%
p(x) $ is free of squares as stated in proposition 1.

More generally, this result completely solves the problem of strong $BF$%
-equivalence for those $\mathbb{Z}\left[ \mathbb{\beta }\right] $ in which
all ideals are invertible in the respective coefficient ring, in particular,
in the quadratic case. In fact, for $n=2$, the semigroup of ideal classes of 
$\mathbb{Z}\left[ \mathbb{\beta }\right] $ has the structure of a Clifford
semigroup, that is, a commutative semigroup $\mathcal{S}$ such that each $%
x\in \mathcal{S}$ is contained in a subgroup of $\mathcal{S}$. The paper 
\textbf{[Za-Zn 94]} contains a detailed description of this situation, as
well as the following result:

\begin{theorem}
(Zanardo and Zannier) For any field $K$ of algebraic numbers of degree $n>2$%
, there is an order $R$ contained in $Z_{K}$ such that the semigroup of
ideal classes of $R$ is not a Clifford semigroup.
\end{theorem}

Proposition 9\ can not be generalized to all automorphisms. However, we
consider the following:

\begin{problem}
Are two $\mathcal{L}$-equivalent automorphisms necessarily $BF$-equivalent?
\end{problem}

Although we are unable to prove this for the moment, we will show that this
is true for a larger class of ideals in the next section.

\subsection{Ideals in monogenic orders}

We start by considering in detail the properties of ideals in the special
case $C(I)=\mathbb{Z}\left[ \theta \right] $ for some algebraic integer $%
\theta $. Orders of this kind are called monogenic.

We need one more definition:

\begin{definition}
Given an ideal $I\vartriangleleft \mathbb{Z}\left[ \mathbb{\beta }\right] $,
the dual ideal of $I$ is 
\begin{equation*}
I^{\ast }=\left\{ z\in K:Tr(zy)\in \mathbb{Z},\forall y\in I\right\} .
\end{equation*}
If $\left( v_{1},\cdots ,v_{n}\right) $ is an integral basis of $I$, the
dual basis is the integral basis $\left( w_{1},\cdots ,w_{n}\right) $ of $%
I^{\ast }$ defined by $Tr(v_{i}w_{j})=\delta _{ij}$.
\end{definition}

This notion may be extended to equivalence classes of ideals. Notice also
that $I^{\ast \ast }=I$.

For the ring $\mathbb{Z}\left[ \mathbb{\beta }\right] $ itself one may
choose the entries of $w_{0}=\left( 1,\mathbb{\beta },\cdots ,\mathbb{\beta }%
^{n-1}\right) $ as the canonical integral basis of $\mathbb{Z}\left[ \mathbb{%
\beta }\right] $. Its dual basis is then $\frac{1}{p^{\prime }(\mathbb{\beta 
})}\left( b_{0},\cdots ,b_{n-1}\right) $ where the $b_{j}$ are another
integral basis of $\mathbb{Z}\left[ \mathbb{\beta }\right] $, and $p^{\prime
}(\mathbb{\beta })$ is just the usual derivative. To be precise, they are
obtained (cf. \textbf{[Fr-Ta 93]}) from the following polynomial identity 
\begin{equation*}
\frac{p(x)}{x-\mathbb{\beta }}=\sum_{j=0}^{n-1}b_{j}x^{j}
\end{equation*}

Denote by $C$ the integer matrix with characteristic polynomial $p(x)$ and
column eigenvector $w_{0}^{t}$. A direct calculation shows that $v_{0}=$ $%
\frac{1}{p^{\prime }(\mathbb{\beta })}\left( b_{0},\cdots ,b_{n-1}\right) $
is exactly the row eigenvector of $C$ normalized such that $v_{0}w_{0}=1$.

In the general case, if $w$ is a column eigenvector of $A$ associated to $%
\mathbb{\beta }$, there exists then a non-singular rational matrix $M$ such
that $w=Mw_{0}$ and directly from the definition we obtain that the dual
basis of $w$ is given by $v=v_{0}M^{-1}$ which is a row eigenvector for $A$,
since $MA=CM$; moreover, $vw=1$ as well. So, $I$ and $I^{\ast }$ are
naturally identified with, resp., the lattice $\mathbb{Z}^{n}$ and $\left( 
\mathbb{T}^{n}\right) ^{\ast }$.

\begin{proposition}
Given $I\vartriangleleft \mathbb{Z}\left[ \mathbb{\beta }\right] $, $%
I^{-1}=p^{\prime }(\beta )I^{\ast }$
\end{proposition}

\begin{proof}
Let $y\in I^{-1}$; then, $\forall z\in I$, we have $Tr(\frac{y}{p^{\prime
}(\beta )}z)\in \mathbb{Z}$, since $yz\in \mathbb{Z}\left[ \mathbb{\beta }%
\right] $ and $\frac{1}{p^{\prime }(\beta )}\in \frac{1}{p^{\prime }(\beta )}%
\mathbb{Z}\left[ \mathbb{\beta }\right] =(\mathbb{Z}\left[ \mathbb{\beta }%
\right] )^{\ast }$. This shows that $I^{-1}\subset p^{\prime }(\beta
)I^{\ast }$. If $u\in I^{\ast }$ then $\forall z\in I$ and $\forall a\in 
\mathbb{Z}\left[ \mathbb{\beta }\right] $ we have $Tr(u(za))\in \mathbb{Z}$
which means that $uI\subset (\mathbb{Z}\left[ \mathbb{\beta }\right] )^{\ast
}=\frac{1}{p^{\prime }(\beta )}\mathbb{Z}\left[ \mathbb{\beta }\right] $ or $%
p^{\prime }(\beta )u\in I^{-1}$ establishing the reverse inclusion
\end{proof}

As a consequence 
\begin{equation*}
(I^{-1})^{-1}=(p^{\prime }(\beta )I^{\ast })^{-1}=\frac{1}{p^{\prime }(\beta
)}(I^{\ast })^{-1}=\frac{1}{p^{\prime }(\beta )}(p^{\prime }(\beta )I^{\ast
}{}^{\ast })=I
\end{equation*}
i.e. any $\mathbb{Z}\left[ \mathbb{\beta }\right] $-ideal is divisorial.

\begin{corollary}
An ideal $I\vartriangleleft \mathbb{Z}\left[ \mathbb{\beta }\right] $ is
invertible iff $C(I)=\mathbb{Z}\left[ \mathbb{\beta }\right] $.
\end{corollary}

\begin{remark}
This confirms also that every order $R$ containing $\mathbb{Z}\left[ \mathbb{%
\beta }\right] $ is the ring of coefficients of some $\mathbb{Z}\left[ 
\mathbb{\beta }\right] $-ideal, namely of $(\mathbb{Z}\left[ \mathbb{\beta }%
\right] :R)$ ($=R^{-1}$ viewing $R$ as a $\mathbb{Z}\left[ \mathbb{\beta }%
\right] $-fractional ideal). In fact, being divisorial, $(R^{-1})^{-1}=R$.
This shows that $C(R^{-1})\subset R$ and the other inclusion is obvious.
\end{remark}

We end this section by describing the place of non-invertible prime ideals
of $\mathbb{Z}\left[ \mathbb{\beta }\right] $.

If $P\vartriangleleft \mathbb{Z}\left[ \mathbb{\beta }\right] $ is a
non-invertible prime ideal, we have $C(P)=P^{-1}$ which is, thus, an order
above $\mathbb{Z}\left[ \mathbb{\beta }\right] $. If $R$ is another order
such that $\mathbb{Z}\left[ \mathbb{\beta }\right] \subset R\subset P^{-1}$,
then, by definition, $P\subset (\mathbb{Z}\left[ \mathbb{\beta }\right]
:R)=R^{-1}$. So, either $R^{-1}=\mathbb{Z}\left[ \mathbb{\beta }\right] $
and $R=\mathbb{Z}\left[ \mathbb{\beta }\right] $ as well, or $R^{-1}=P$ and,
by divisoriality, $R=P^{-1}$. The same reasoning shows the converse: if $R$
is an order immediately above $\mathbb{Z}\left[ \mathbb{\beta }\right] $,
then $(\mathbb{Z}\left[ \mathbb{\beta }\right] :R)=R^{-1}$ is a prime ideal
of $\mathbb{Z}\left[ \mathbb{\beta }\right] $.

\subsection{$BF$-equivalence: a partial generalization}

By the reasoning preceding the last proposition, in the monogenic case the
matrix representing multiplication by $\beta $, with respect to an
appropriate basis of $I^{-1}$, will be $A^{t}$, which implies that they are $%
BF$-equivalent. So, just by transposing both sides of the direct sum
conjugacy defined in the preceding section, we get the

\begin{proposition}
If two automorphisms $A$ and $A^{\prime }$ with associated ideals $I$ and $J$
are $\mathcal{L}$-equivalent and at least one ideal in each of the pairs $(I$
, $(\mathbb{Z}\left[ \mathbb{\beta }\right] :I))$ and $(J$ , $(\mathbb{Z}%
\left[ \mathbb{\beta }\right] :J))$ is invertible, then $A$ and $A^{\prime }$
are $BF$-equivalent.
\end{proposition}

We'll now show that this effectively enlarges the class of $BF$-equivalent
automorphisms. We need another bit of information on the role of the ring of
coefficients $C(I)=R$ for $I\vartriangleleft \mathbb{Z}\left[ \mathbb{\beta }%
\right] $.

\emph{All inverses in the next sentence must be considered with respect to }$%
Z\left[ \mathbb{\beta }\right] $. Lemma 6, applied to this case,\ says that $%
(II^{-1})^{-1}=C(I^{-1})=C((I^{-1})^{-1})=C(I)=R,$ because $I$ is divisorial
as a $\mathbb{Z}\left[ \mathbb{\beta }\right] $-ideal. By the same reason, $%
II^{-1}=((II^{-1})^{-1})^{-1}=R^{-1}$, i.e.

\begin{lemma}
Given $I\vartriangleleft \mathbb{Z}\left[ \mathbb{\beta }\right] $ and $(%
\mathbb{Z}\left[ \mathbb{\beta }\right] :I)$ its inverse, $I(\mathbb{Z}\left[
\mathbb{\beta }\right] :I)=(\mathbb{Z}\left[ \mathbb{\beta }\right] :R)$,
where $R=C(I)$.
\end{lemma}

If some $I$ with $C(I)=R$ is not invertible, then $(\mathbb{Z}\left[ \mathbb{%
\beta }\right] :R)$ is also a non-invertible ideal with the same coefficient
ring. On the other hand, if $(\mathbb{Z}\left[ \mathbb{\beta }\right] :R)$
is not invertible, for every ideal $I$, with $C(I)=R$, at least one of $I$
and $(\mathbb{Z}\left[ \mathbb{\beta }\right] :I)$ is also not invertible.

So, given an order $R$ above $\mathbb{Z}\left[ \mathbb{\beta }\right] $,
there are two possible cases: either all ideals with $R$ as coefficient ring
are invertible (this will be the case if $R$ is monogenic) and we know that
their associated automorphisms are all strongly $BF$-equivalent; or there
are non-invertible ideals with coefficient ring $R$; in this case, for each
invertible ideal $I$, $(\mathbb{Z}\left[ \mathbb{\beta }\right] :I)$ is not
invertible.

In case the converse statement is true (this would be: if $I$ with $C(I)=R$
is not invertible then $(\mathbb{Z}\left[ \mathbb{\beta }\right] :I)$ is
invertible), proposition 11 would solve the problem. However, it remains the
possibility of both $I$ and $(\mathbb{Z}\left[ \mathbb{\beta }\right] :I)$
being non-invertible.

\section{BF-groups and flow equivalence}

We end this paper by returning, in some way, to the original role of
BF-groups, i.e. its connection to flow equivalence. We start with the
necessary definitions.

\begin{definition}
Let $f:X\rightarrow X$ be a homeomorphism of a topological space. The
suspension of $X$ determined by $f$ is the quotient space $Y=(X\times 
\mathbb{R})/\sim $ where the equivalence relation is given by $\left(
x,t+1\right) \sim \left( f(x),t\right) $.
\end{definition}

The projection from $X\times \mathbb{R}$ to $Y$ defines a flow in the later
space and each flow line contains exactly one orbit of $f$.

\begin{definition}
Two homeomorphisms $f_{1}:X_{1}\rightarrow X_{1}$ and $f_{2}:X_{2}%
\rightarrow X_{2}$ are flow-equivalent if there exists an homeomorphism $%
\varphi :Y_{1}\rightarrow Y_{2}$ conjugating the two induced flows and
preserving its direction.
\end{definition}

The Bowen-Franks group $BF_{1}$ was introduced by the two authors as a
flow-equivalence invariant in the case where the spaces $X_{i}$ are
subshifts of finite type. Later, the second author showed that it is an
almost complete invariant:

\begin{theorem}
Let $(X_{A},\sigma )$ and $(X_{B},\sigma )$ be two irreducible subshifts of
finite type with positive entropy, where $A,B$ are the transition matrices.
Then the shifts acting in the two spaces are flow equivalent iff $\det
(Id-A)=\det (Id-B)$ and $BF_{1}(A)=BF_{1}(B)$.
\end{theorem}

We want to show how the Bowen-Franks group appears naturally as a
flow-equivalence invariant in the case of torus automorphisms. Let $A$ be a $%
\mathbb{T}^{n}$-automorphism and $Y$ the suspension space it determines.

The fundamental group of $\mathbb{T}^{n}\times \mathbb{R}$, a free abelian
group of rank $n$, is embedded into the fundamental group of $Y$. Let $%
x_{1},\ldots ,x_{n}$ be a basis of this subgroup of $\pi _{1}(Y)$. A new
generator $x_{0}$ must be added, namely the one represented by the closed
orbit through $0$, i.e. the loop$\left\{ \left( 0,t\right) :t\in \mathbb{R}%
\right\} $.

The presentation of $\pi _{1}(Y)$ is completed with the relations imposed by
the equivalence $\left( x,t\right) \sim \left( Ax,t-1\right) $. In order to
describe these in a simple way, we denote $x_{1}^{m_{1}}\cdots x_{n}^{m_{n}}$
by $X^{m}$ where $m=\left( m_{1},\ldots ,m_{n}\right) $.

We leave to the reader the verification of the following

\begin{proposition}
$\pi _{1}(Y)=\{x_{0},x_{1},\ldots ,x_{n}:x_{i}x_{j}=x_{j}x_{i},\forall 1\leq
i,j\leq n;x_{0}X^{m}=X^{mA}x_{0},\forall m\in \mathbb{Z}^{n}\}$.
\end{proposition}

The abelianization of $\pi _{1}(Y)$ gives

\begin{corollary}
$H_{1}(Y)=\mathbb{Z\oplus Z}^{n}/\mathbb{Z}^{n}(A-I)=\mathbb{Z\oplus }%
BF_{1}(A)$.
\end{corollary}

As an homeomorphism between the suspension spaces determined by different $%
\mathbb{T}^{n}$-automorphisms induces an isomorphism between homology
groups, it's obvious that the Bowen-Franks group $BF_{1}$ is an invariant of
flow-equivalence.

\section*{Acknowledgments}

We would like to thank FCT (Portugal) for having in part supported this work
through the Research Units Pluriannual Funding Program.

\end{document}